\documentclass[10pt]{amsart}
\usepackage{amssymb}

\usepackage{amscd}
\usepackage{amsthm}
\usepackage[dvips]{graphicx}

\usepackage[all]{xy}

\newtheorem{Thm}{Theorem}[section]
\newtheorem{Lem}[Thm]{Lemma}

\newtheorem{Cor}[Thm]{Corollary}
\newtheorem{Prop}[Thm]{Proposition}

\newtheorem{Rem1}[Thm]{Remark}

\title{Auslander-Reiten conjecture for symmetric algebras of polynomial growth}
\author{ Guodong Zhou and Alexander Zimmermann}

\address{Guodong Zhou
\newline Institut f\"ur Mathematik,
\newline Universit\"at Paderborn,
\newline 33095 Paderborn,
\newline Germany}
\email{gzhou@math.uni-paderborn.de}
\address{Alexander Zimmermann
\newline Universit\'e de Picardie,
\newline D\'epartement de Math\'ematiques et LAMFA (UMR 6140 du CNRS),
\newline 33 rue St Leu,
\newline F-80039 Amiens Cedex 1,
\newline France}
\email{alexander.zimmermann@u-picardie.fr}
\date{May 3, 2010}

\keywords{Auslander--Reiten conjecture; Derived equivalence;  Reynolds
ideal;  Self-injective algebras of polynomial growth; Stable
equivalence; Stable equivalence of Morita type}

\subjclass[2000]{16G10, 16E40, 20C20}

\newenvironment{Proof}[1][Proof]{\begin{trivlist}
\item[\hskip \labelsep {\bfseries #1}]}{\flushright
$\Box$\end{trivlist}}

\newenvironment{Rem}{\begin{Rem1}\rm}{\end{Rem1}}

\newcommand{\lra}{\longrightarrow}

\newcommand{\ra}{\rightarrow}
\newcommand{\sdp}{\times\kern-.2em\vrule height1.1ex depth-.05ex}
\newcommand{\epi}{\lra \kern-.8em\ra}

\newcommand{\N}{{\mathbb N}}

\newcommand{\ul}{\underline}

 \normalbaselines
\begin{document}

\begin{abstract} This paper studies  self-injective algebras of
polynomial growth.  We prove that the derived equivalence
classification of weakly symmetric algebras of domestic type
coincides with the classification up to stable equivalences (of
Morita type).  As for weakly symmetric non-domestic algebras of
polynomial growth,  up to some scalar problems, the derived
equivalence classification coincides with the classification up to
stable equivalences of Morita type.  As a consequence, we get the
validity of  the Auslander-Reiten conjecture for stable
equivalences of Morita type between weakly symmetric algebras of
polynomial growth.
\end{abstract}

\maketitle

\section*{Introduction}

Drozd \cite{Drozd}
showed that finite dimensional algebras $A$ over an algebraically closed
field $K$ are either of finite representation type, or of tame representation type
or of wild representation type.
Here, an algebra is of finite representation type if there are only a finite number
of isomorphism classes of indecomposable modules. An algebra is of
wild representation type if all module categories are full subcategories
of the module category of the algebra. An algebra is of tame representation
type if for every positive integer $d$ there are $A-K[X]$-bimodules $M_i(d)$
for $i\in\{1,\dots,n_d\}$ so that $M_i(d)$ is free as $K[X]$-module and so that
all but finite number of isomorphism classes of indecomposable
$d$-dimensional $A$-modules are isomorphic to
$M_i(d)\otimes_{K[X]}K[X]/(X-\lambda)$ for $\lambda\in K$ and $i\in\{1,\dots,n_d\}$.
If $n_d$ is minimal with respect to $d$ satisfying this property,
$A$ is of polynomial growth if there is an integer $m$ so that
$\lim_{d\ra\infty}\frac{n_d}{d^m}=0$, and $A$ is domestic if there is an integer
$m$ so that $n_d\leq d\cdot m$ for all $d$.

Bocian, Holm and Skowro\'nski
\cite{BocianHolmSkowronski2004,BocianHolmSkowronski2007,HolmSkowronski2006}
classified all weakly symmetric
algebras of domestic type up to derived equivalence and Bia\l kowski,
Holm and Skowro\'{n}ski \cite{BialkowskiHolmSkowronski2003a,
BialkowskiHolmSkowronski2003b,HolmSkowronski2009}
did the same thing for all weakly symmetric
 non-domestic
algebras of polynomial growth. The result is a finite number of
families of algebras given by quivers and relations, depending on
certain parameters. In the case  of domestic algebras the
classification is complete, whereas in case of weakly symmetric
algebras of polynomial growth it is not known if certain choices
of parameters may lead to the same derived equivalence class.
This question if certain parameters lead to derived equivalent algebras
will be called the scalar problem in the sequel.

The main result of this paper is that the derived equivalence
classification of weakly symmetric algebras of polynomial growth
coincides (up to some scalar problems) with the
classification up to stable equivalences of
Morita type as given in
\cite{BocianHolmSkowronski2004,BocianHolmSkowronski2007,
HolmSkowronski2006} and \cite{BialkowskiHolmSkowronski2003a,
BialkowskiHolmSkowronski2003b,HolmSkowronski2009}. As a
consequence, the Auslander-Reiten conjecture holds for a stable
equivalence of Morita type between two weakly symmetric algebras of
polynomial growth.

For weakly symmetric algebras of domestic type we can show more. We prove that
the derived equivalence classification coincides with the classification up
to stable equivalences (no scalar problem occurs in this case) and that  the Auslander-Reiten conjecture holds
for stable equivalences between weakly symmetric algebras of domestic types.

As mentioned above the classification
of symmetric algebras of polynomial growth up to derived equivalence
leaves open if for certain choices of parameters the
algebras are derived equivalent or not.
The same questions concerning
the same choices of parameters remain open also for the
classification up to stable equivalences of Morita type.

Moreover, using our result \cite{ZhouZimmermannTameBlocks}, we show that
for tame symmetric algebras with periodic modules,
the derived equivalence classification coincides with the classification up
to stable equivalence of Morita type, up to some scalar problem, and that  the Auslander-Reiten conjecture holds
for stable equivalences between tame symmetric algebras with periodic modules.

The paper is organised as follows. In Section~\ref{stableinvsect}
we recall basic definitions and recall
the stable invariants we use in the sequel. In
Section~\ref{weaklysymmdomesticsection} we prove the
main result Theorem~\ref{domesticmaintheorem}
for weakly symmetric algebras of domestic representation type and in
Section~\ref{polysection} we show the main result
Theorem~\ref{NonDomesticClassification}
for stable equivalences of Morita type for weakly
symmetric algebras of polynomial
growth.

As soon as we use statements on representation type we shall always assume that
$K$ is an algebraically closed field.

\section{Stable categories; definitions, notations and invariants}

\label{stableinvsect}

\subsection{The definitions}

Let $K$ be a field and let $A$ be a finite dimensional
$K$-algebra. The stable category $A-\ul{mod}$ has
as objects all finite dimensional $A$-modules and  morphisms
from na $A$-module $X$ to another $A$-module $Y$ are the
equivalence classes
of $A$-linear homomorphisms $Hom_A(X,Y)$, where
two morphisms of $A$-modules are declared to be equivalent if
their difference factors through a projective $A$-module.

Two algebras are called {\em stably equivalent}
if $A-\ul{mod}\simeq B-\ul{mod}$ as additive categories.
Let $M$ be an $A-B$-bimodule and let $N$ be a
$B-A$-bimodule. Following Brou\'e the couple $(M,N)$
defines a {\em stable equivalence of Morita type} if
\begin{itemize}
\item
$M$ is projective as $A$-module and is projective as $B$-module
\item
$N$ is projective as $B$-module and is projective as $A$-module
\item
$M\otimes_BN\simeq A\oplus P$ as $A-A$-bimodules, where $P$ is a projective
$A-A$-bimodule.
\item
$N\otimes_AM\simeq B\oplus Q$ as $B-B$-bimodules, where $Q$ is a projective
$B-B$-bimodule.
\end{itemize}
Of course, if $(M,N)$ defines a stable equivalence of Morita type,
then $M\otimes_B-:B-\ul{mod}\simeq A-\ul{mod}$ is an equivalence. However,
it is known that there are algebras which are stably equivalent but for which there
is no stable equivalence of Morita type.

\subsection{Invariants under stable equivalences}

Auslander and Reiten conjectured that if $A$ and $B$ are stably equivalent,
then the number of isomorphism classes of
 non-projective simple $A$-modules is equal to the number of
isomorphism classes of  non-projective simple $B$-modules. The conjecture
is open in general, but we shall use in this article the following positive
result \cite{Pogorzaly} of
Pogorza\l y: If $A$ is stably equivalent to $B$ and if $A$ is selfinjective special biserial that is not Nakayama,
then $B$ is self-injective special biserial as well and the Auslander-Reiten
conjecture holds in this case.

Krause and Zwara \cite{KrauseZwara} showed that stable equivalences preserve the representation type.
Suppose $A$ and $B$ are two finite dimensional stably equivalent
$K$-algebras. If $A$ is of domestic representation
type (respectively of polynomial growth, respectively of tame
representation type), then so is $B$.

Reiten \cite{Reiten} show that being selfinjective is almost
invariant under stable equivalences. More precisely, let $A$ be a selfinjective
finite dimensional $K$-algebra which is stably equivalent
to a  finite dimensional indecomposable $K$-algebra $B$.
Then $B$ is either selfinjective or a radical squared zero Nakayama algebra.

Finally, if $A$ and $B$ are stably equivalent finite dimensional
algebras with Loewy length at least three,
then the stable Auslander-Reiten quivers are isomorphic as translation quivers.

\subsection{Invariants under stable equivalence of Morita type}

Much more is known for stable equivalences of Morita type. We shall cite only
those properties used in the sequel.

Xi has shown \cite{Xi2008} that if two finite dimensional $K$-algebras
$A$ and $B$ are stably
equivalent of Morita type, then the absolute value of the Cartan determinants
are equal. Even better, the elementary divisors, including their multiplies,
of the Cartan matrices different from $\pm 1$ coincide for $A$ and for $B$.

Liu showed \cite{Liu2008}
that if $A$ is an indecomposable finite dimensional $K$-algebra,
and $B$ is a finite dimensional $K$-algebra stably
equivalent of Morita type to $A$, then $B$ is indecomposable as well.

If $K$ is an algebraically closed field of characteristic $p>0$ and $A$ a finite dimensional
$K$-algebra. Define $[A,A]$ the $K$-space generated by all elements $ab-ba\in A$
for all $a,b\in A$. Then, $$\{x\in A\;|\;x^{p^n}\in[A,A]\}=:T_n(A)$$
is a $K$-subspace of $A$ containing $[A,A]$.
Liu and the authors have shown \cite{LZZ}
that
$$dim_K(T_n(A)/[A,A])=dim_K(T_n(B)/[B,B])$$
for all $n\in\N$. If $A$ and $B$ are symmetric and stably equivalent of
Morita type, then taking the orthogonal space with respect to the symmetrising
form, observing that $[A,A]^\perp=Z(A)$, Liu, K\"onig and the first author have
shown \cite{KLZ} that $$Z(A)/T_n(A)^\perp\simeq Z(B)/T_n(B)^\perp$$
as rings.

\section{Weakly symmetric  algebras of domestic representation type }

\label{weaklysymmdomesticsection}

Selfinjective algebras of domestic type split
into two subclasses: standard algebras and non-standard ones.

A selfinjective algebra of tame representation type is called {\em standard} if its
basic algebra admit simply connected Galois coverings. A
selfinjective algebra of tame representation type is called {\em nonstandard}
if it is not standard.

\subsection{Standard domestic weakly symmetric algebras}

\begin{Thm}\cite[Theorem 1]{BocianHolmSkowronski2004}
\label{DomesticStandardSingular}
For an algebra $A$ the following
statements are equivalent:
\begin{itemize}
\item[(1)] $A$ is representation-infinite domestic weakly
symmetric algebras having simply connected Galois coverings and
the Cartan matrix $C_A$ is singular.

\item[(2)] $A$ is derived equivalent to the trivial extension
$T(C)$ of a canonical algebra $C$ of Euclidean type.

\item[(3)] $A$ is stably equivalent to the trivial extension
$T(C)$ of a canonical algebra $C$ of Euclidean type.
\end{itemize}
Moreover, the trivial extensions $T(C)$ and $T(C')$ of two
canonical algebras $C$ and $C'$ of Euclidean type are derived
equivalent (respectively, stably equivalent) if and only if the
algebras $C$ and $C'$ are isomorphic.
\end{Thm}

From this theorem, we know that   a standard weakly symmetric
algebras of domestic representation type with {\em singular Cartan
matrice} is symmetric.

\begin{Cor} The class of indecomposable standard weakly
symmetric  algebras of domestic representation type with singular
Cartan matrices is closed under stable equivalences.

\end{Cor}

\begin{Proof}
 In fact, let $A$ be an indecomposable algebra stably 
 equivalent to an indecomposable standard weakly
symmetric  algebras of domestic representation type. 
Then by the preceding theorem, $A$ is stably equivalent to
the trivial extension
$T(C)$ of a canonical algebra $C$ of Euclidean type. 
Again by the preceding theorem, $A$ is standard weakly
symmetric  algebras of domestic representation type.
\end{Proof}

For standard weakly symmetric algebras of domestic
representation type with {\em nonsingular Cartan matrices,}  we have the
following derived norm forms.
$$\unitlength0.6cm
\begin{picture}(11,6)
 \put(9.4,3){\circle{2.0}}
\put(6.8,3){\circle{2.0}} \put(7.85,2.9){\vector(0,1){0.3}}
\put(8,2.9){$\bullet$} \put(8.35,2.9){\vector(0,1){0.3}}
\put(4.8,2.9){$\beta$} \put(11,2.9){$\alpha$}

 \put(0,1.9){$A(\lambda)$}
\put(0,1){$\lambda\in K\backslash \{0\}$}

\put(4.8,0.5){$\alpha^2=0, \beta^2=0,
\alpha\beta=\lambda\beta\alpha$}
\end{picture}$$

$$\unitlength0.5cm
\begin{picture}(30,15)
 \put(18.1, 11.1){\vector(1,1){1.8}}

\put(16.5, 11.5){$\beta_1$}

 \put(18, 11){\vector(-1,1){1.8}}

 \put(14.5, 13){$\beta_2$}

\put(16, 13){\vector(-2,1){1.8}}

\put(13, 13.2){$\beta_3$}

\put(14, 14){\vector(-1,0){1.8}}

\put(11, 13.0){$\beta_4$}

\put(12, 14){\vector(-2,-1){1.8}}

\put(10.7, 8.8){$\beta_{q-3}$}

\put(10, 9){\vector(2,-1){1.8}}

\put(12.2, 8.5){$\beta_{q-2}$}

\put(12, 8){\vector(1,0){1.8}}

\put(14, 8.8){$\beta_{q-1}$}

\put(14, 8){\vector(2,1){1.8}}

\put(16.4, 10.3){$\beta_{q}$}

\put(16, 9){\vector(1,1){1.8}}

\put(19.2, 11.5){$\alpha_{1}$}

\put(20, 9){\vector(-1,1){1.8}}

\put(21, 13){$\alpha_{2}$}

\put(22, 8){\vector(-2,1){1.8}}

\put(22.5, 13.5){$\alpha_{3}$}

\put(24, 8){\vector(-1,0){1.8}}

\put(24.5, 13){$\alpha_{4}$}

\put(26, 9){\vector(-2,-1){1.8}}

\put(18.8, 10.3){$\alpha_{p}$}

\put(20.7, 8.8){$\alpha_{p-1}$}

\put(22.4, 8.5){$\alpha_{p-2}$}

\put(24.1, 8.8){$\alpha_{p-3}$}

\put(20, 13){\vector(2,1){1.8}}

\put(22, 14){\vector(1,0){1.8}}

\put(24, 14){\vector(2,-1){1.8}}

\multiput(9, 11)(0.1,0.2){10}{\circle*{0.01}}

\multiput(9, 11)(0.1,-0.2){10}{\circle*{0.01}}

\multiput(27, 11)(-0.1,0.2){10}{\circle*{0.01}}

\multiput(27, 11)(-0.1,-0.2){10}{\circle*{0.01}}

 \put(1, 12){$A(p, q)$}

 \put(1, 10.5){$1\leq p\leq q$}
 \put(1, 9){$p+q\geq 3$}

 \put(9, 6){$\alpha_1\alpha_2\cdots \alpha_p\beta_1\beta_2\cdots
 \beta_q=\beta_1\beta_2\cdots \beta_q\alpha_1\alpha_2\cdots \alpha_p$}

\put(9, 4.5){$ \alpha_p\alpha_1=0, \beta_q\beta_1=0$}

\put(9, 3){$\alpha_i\alpha_{i+1}\cdots \alpha_p\beta_1\cdots
\beta_q  \alpha_1 \cdots \alpha_{i-1}\alpha_i=0, 2\leq i\leq p$}

\put(9, 1.5){$\beta_j\beta_{j+1}\cdots \beta_q\alpha_1\cdots
\alpha_p \beta_1 \cdots \beta_{j-1}\beta_j=0, 2\leq j\leq q$}

\end{picture}$$

$$\unitlength0.5cm
\begin{picture}(30,11)

\multiput(9, 7)(0.1,0.2){10}{\circle*{0.01}}

\multiput(9, 7)(0.1,-0.2){10}{\circle*{0.01}}

\put(16.5, 7.5){$\beta_1$}

 \put(18, 7){\vector(-1,1){1.8}}

 \put(14.5, 9){$\beta_2$}

\put(16, 9){\vector(-2,1){1.8}}

\put(13, 9.2){$\beta_3$}

\put(14, 10){\vector(-1,0){1.8}}

\put(11, 9.0){$\beta_4$}

\put(12, 10){\vector(-2,-1){1.8}}

\put(10.7, 4.8){$\beta_{n-3}$}

\put(10, 5){\vector(2,-1){1.8}}

\put(12.2, 4.5){$\beta_{n-2}$}

\put(12, 4){\vector(1,0){1.8}}

\put(14, 4.8){$\beta_{n-1}$}

\put(14, 4){\vector(2,1){1.8}}

\put(16.4, 6.3){$\beta_{n}$}

\put(16, 5){\vector(1,1){1.8}}

\put(19.2, 7){\circle{2}}

\put(19, 8){\vector(-1,0){0.01}}

\put(20.5, 7){$\alpha$}

 \put(1, 8){$\Lambda(n)$}

 \put(1, 6.5){$n\geq 2$}

 \put(9, 2){$\alpha^2=(\beta_1\beta_2\cdots\beta_n)^2, \alpha\beta_1=0, \beta_n\alpha=0$}

\put(9, 0.5){$\beta_j\beta_{j+1}\cdots \beta_n\beta_1\cdots
\beta_n \beta_1 \cdots \beta_{j-1}\beta_j=0, 2\leq j\leq n$}

\end{picture}$$

$$\unitlength0.5cm
\begin{picture}(30,15)

\multiput(9, 11)(0.1,0.2){10}{\circle*{0.01}}

\multiput(9, 11)(0.1,-0.2){10}{\circle*{0.01}}

\put(16.5, 11.5){$\beta_1$}

 \put(18, 11){\vector(-1,1){1.8}}

 \put(14.5, 13){$\beta_2$}

\put(16, 13){\vector(-2,1){1.8}}

\put(13, 13.2){$\beta_3$}

\put(14, 14){\vector(-1,0){1.8}}

\put(11, 13.0){$\beta_4$}

\put(12, 14){\vector(-2,-1){1.8}}

\put(10.7, 8.8){$\beta_{n-3}$}

\put(10, 9){\vector(2,-1){1.8}}

\put(12.2, 8.5){$\beta_{n-2}$}

\put(12, 8){\vector(1,0){1.8}}

\put(14, 8.8){$\beta_{n-1}$}

\put(14, 8){\vector(2,1){1.8}}

\put(16.4, 10.3){$\beta_{n}$}

\put(16, 9){\vector(1,1){1.8}}

\put(18.2, 11.3){\vector(1,2){1.2}}

\put(19.8, 13.4){\vector(-1,-2){1.2}}

\put(18.4, 10.9){\vector(1,-2){1.2}}

\put(19.2, 8.4){\vector(-1,2){1.2}}

\put(19.4, 9.7){$\gamma_1$}

\put(19.4, 11.7){$\alpha_2$}

\put(18, 9){$\gamma_2$}

\put(18, 12.5){$\alpha_1$}

 \put(1, 12){$\Gamma(n)$}

 \put(1, 10.5){$n\geq 1$}

 \put(9, 6){$\alpha_1\alpha_2=(\beta_1\beta_2\cdots\beta_n)^2=\gamma_1\gamma_2,$}

\put(9, 4.5){$\alpha_2\beta_1=0, \gamma_2\beta_1=0,
\beta_n\alpha_1=0$}

\put(9, 3){$\beta_n\gamma_1=0, \alpha_2\gamma_1=0,
\gamma_2\alpha_1=0$}

\put(9, 1.5){$\beta_j\beta_{j+1}\cdots \beta_n\beta_1\cdots
\beta_n \beta_1 \cdots \beta_{j-1}\beta_j=0, 2\leq j\leq n$}

\end{picture}$$

\begin{Thm}\cite[Theorem
2]{BocianHolmSkowronski2004}\label{DomesticStandardNonSingular}
For a domestic standard self-injective algebra $A$,    the
following statements are equivalent:
\begin{itemize}
\item[(1)]$A$ is weakly symmetric and the Cartan matrix $C_A$ is
nonsingular.

\item[(2)] $A$ is derived equivalent to an algebra of the form
$A(\lambda), A(p, q), \Lambda(n),  \Gamma(n)$.

\item[(3)] $A$ is stably equivalent to an algebra of the form
$A(\lambda), A(p, q), \Lambda(n), \Gamma(n)$.

\end{itemize}

Moreover, two algebras of the forms $A(\lambda), A(p, q),
\Lambda(n)$ or $\Gamma(n)$ are derived equivalent (respectively,
stably equivalent) if and only if they are isomorphic.
\end{Thm}
All these algebras except
$$A(\lambda)=k\langle X, Y\rangle/(X^2, Y^2, XY-\lambda YX), \lambda\not\in\{ 0, 1\}$$
are symmetric. Remark that except $\Gamma(n)$, all algebras are
special biserial algebras.

Since the class of indecomposable standard weakly
symmetric  algebras of domestic representation type is closed under stable equivalences, we have the following

\begin{Cor} \label{standardweaklysymmetricstablyequiv}
Two  standard weakly symmetric algebras of
domestic representation type are stably equivalent if and only
if they are derived equivalent.
\end{Cor}

\subsection{Non standard domestic weakly symmetric algebras}

Nonstandard self-injective algebras of domestic
representation type are classified up to derived and stable
equivalences in \cite{BocianHolmSkowronski2007}.

\begin{Thm}\cite[Theorem 1]{BocianHolmSkowronski2007}\label{DomesticNonStandard}
Any nonstandard representation-infinite
selfinjective algebra of domestic type is derived equivalent
(resp. stably equivalent) to an algebra $\Omega(n)$ with $n\geq
1$. Moreover, two algebras $\Omega(n)$ and $\Omega(m)$ are derived
equivalent (respectively, stably equivalent) if and only if $m =
n$.\end{Thm} The quiver with relations of $\Omega(n)$ is as
follows.

$$\unitlength0.5cm
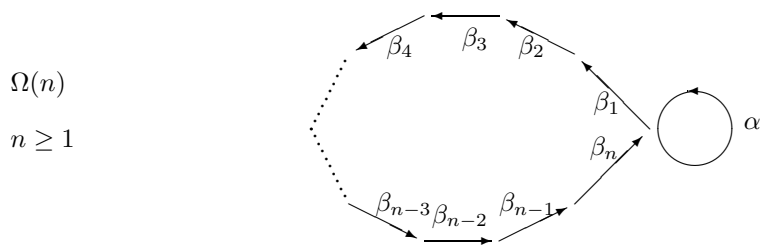
\begin{picture}(30,11)

\multiput(9, 7)(0.1,0.2){10}{\circle*{0.01}}

\multiput(9, 7)(0.1,-0.2){10}{\circle*{0.01}}

\put(16.5, 7.5){$\beta_1$}

 \put(18, 7){\vector(-1,1){1.8}}

 \put(14.5, 9){$\beta_2$}

\put(16, 9){\vector(-2,1){1.8}}

\put(13, 9.2){$\beta_3$}

\put(14, 10){\vector(-1,0){1.8}}

\put(11, 9.0){$\beta_4$}

\put(12, 10){\vector(-2,-1){1.8}}

\put(10.7, 4.8){$\beta_{n-3}$}

\put(10, 5){\vector(2,-1){1.8}}

\put(12.2, 4.5){$\beta_{n-2}$}

\put(12, 4){\vector(1,0){1.8}}

\put(14, 4.8){$\beta_{n-1}$}

\put(14, 4){\vector(2,1){1.8}}

\put(16.4, 6.3){$\beta_{n}$}

\put(16, 5){\vector(1,1){1.8}}

\put(19.2, 7){\circle{2}}

\put(19, 8){\vector(-1,0){0.01}}

\put(20.5, 7){$\alpha$}

 \put(1, 8){$\Omega(n)$}

 \put(1, 6.5){$n\geq 1$}

 \put(3, 2){$\alpha^2=\alpha\beta_1\beta_2\cdots\beta_n,
 \alpha\beta_1\beta_2\cdots\beta_n+\beta_1\beta_2\cdots\beta_n\alpha=0, $}

\put(3, 0.5){$\beta_n\beta_{1}=0, \beta_j\beta_{j+1}\cdots
\beta_n\beta_1\cdots \beta_n \beta_1 \cdots \beta_{j-1}\beta_j=0,
2\leq j\leq n$}

\end{picture}$$

Remark that the algebra $\Omega(n)$ is  always weakly symmetric,
but it is symmetric only when the characteristic of the base field
is $2$. Note that $\Omega(n)$ is not special biserial.

\subsection{Domestic standard type   versus domestic non standard type}

\begin{Thm}\cite[Theorem 1.1]{HolmSkowronski2006}\label{SymmetricDomesticStandardVsNonStandard} A
standard symmetric algebra of domestic representation type cannot
be derived equivalent to a non-standard symmetric one.

\end{Thm}
 In course of the proof of the above theorem,  one needs to compare the algebras
$\Omega(n)$ with $A(1, n)$ in case the characteristic of $K$ is
$2$. It is proved in \cite{HolmSkowronski2006} that the dimensions
of   the centre modulo the first K\"{u}lshammer ideal,  for
$\Omega(n)$ and  $A(1, n)$, are different, but this dimension is
in fact invariant under stable equivalences of Morita type (\cite{LZZ}). So  a
standard symmetric algebra of domestic representation type cannot
be stably equivalent of Morita type to a non-standard symmetric
one. In fact, we can even prove more.

\begin{Lem}\label{WeaklySymmetricDomesticStandardVsNonStandard}
 A standard weakly symmetric algebra of domestic type
cannot be stably equivalent to a nonstandard one.
\end{Lem}

\begin{Proof}
Stably equivalent algebras of Loewy length at least $3$ have
isomorphic stable Auslander-Reiten quivers (\cite[Corollary X.
1.9]{AuslanderReiten}). Notice that the stable Auslander-Reiten
quiver of $A(\lambda)$ for $\lambda\neq 0$ consists of an
Euclidean component of type $\mathbb{Z}\tilde{A}_{1}$ and a
$\mathbb{P}_1(K)$-family of homogenous tubes. Since all algebras
in the above list of Theorems~\ref{DomesticStandardSingular},
\ref{DomesticStandardNonSingular}, \ref{DomesticNonStandard} are
of Loewy length at least $3$,
 by comparing the shapes of the stable
Auslander-Reiten quivers  as in \cite[Section
4]{HolmSkowronski2006}, we know that a standard weakly symmetric
algebra of domestic type $A$ in the list of
Theorems~\ref{DomesticStandardSingular},
\ref{DomesticStandardNonSingular} is stably equivalent to a
nonstandard one $B=\Omega(n)$ only when $A=A(1, n)$ and
$B=\Omega(n)$ for $n\geq 1$. However, notice that $A(1, n)$ is
special biserial and $\Omega(n)$ is not. By the result of
Pogorza\l y \cite[Theorem 7.3]{Pogorzaly}, they can never be
stably equivalent.
\end{Proof}

\subsection{Stable equivalence classification of domestic weakly symmetric algebras}

 \begin{Prop} \label{ClassificationSymmetricDomestic}
 \begin{enumerate}
 \item
 Two symmetric  algebras of domestic representation type   are
derived equivalent
 if and only if they are stably equivalent of Morita type if and only
 if they are stably equivalent.

 \item The class of    symmetric
algebras  of domestic representation type  is closed under stable
equivalences, hence is closed under stably equivalences of Morita type and
derived equivalences.
\end{enumerate}
\end{Prop}

\begin{Proof}
Corollary~\ref{standardweaklysymmetricstablyequiv} shows that for
standard weakly symmetric algebras the notions of derived
equivalence and of stable equivalence are the same.
Lemma~\ref{WeaklySymmetricDomesticStandardVsNonStandard} shows
that standard weakly symmetric algebras cannot be stably
equivalent to non standard weakly symmetric algebras.
Theorem~\ref{DomesticNonStandard} shows that for non standard
domestic weakly symmetric algebras derived and stable equivalence
are the same notion. This proves (1).

\medskip

Now let $A$ be an indecomposable algebra stably equivalent to an
indecomposable  symmetric algebra $B$ of domestic representation
type in the list of algebras in
Theorems~\ref{DomesticStandardSingular},
\ref{DomesticStandardNonSingular} and
 \ref{DomesticNonStandard}. Then by a result of Idun Reiten
\cite[Theorem 2.6]{Reiten} (see also
\cite[Theorem]{LiuPaquette}), $A$ itself is also self-injective or
$A$ is isomorphic to a radical square zero algebra  whose quiver
is of type $A$ with linear orientation. But the latter algebra has
finite representation type, and thus cannot be stably equivalent
to a representation-infinite algebra. We infer that $A$ is
self-injective and now  by \cite[Proposition
5.2]{BocianHolmSkowronski2004}, $A$ is weakly symmetric. By
\cite[Corollary 2 and the discussions afterwards]{KrauseZwara},
$A$ is of domestic type. We thus proved that $A$ is weakly
symmetric of domestic type.

To finish the proof of the second assertion, one needs to show that
$A$ is symmetric and thus to  exclude the cases where $A$ is
stably equivalent to $A(\lambda)$ with $\lambda\in K\backslash\{0,
1\}$ or to $\Omega(n)$ with $n\geq 1$ in case of $char K\neq 2$,
since these are the only cases where non symmetric algebras may
occur.

If $A$ is stably equivalent to $A(\lambda)$ with $\lambda\in
K\backslash\{0, 1\}$, then by  \cite[Theorem 7.3]{Pogorzaly}, $B$
is special biserial. Since the Auslander-Reiten conjecture is
proved for self-injective special biserial algebras (\cite[Theorem
0.1]{Pogorzaly}), $B$ is also local and is thus necessary
isomorphic to $A(1)$. However,  $A(1)$  is never stably
equivalent to $A(\lambda)$ with $\lambda\neq 1$ by an unpublished
result of Rickard.

If $char(K)\neq 2$ and   $A$ is stably equivalent to $A=\Omega(n)$
with $n\geq 1$,  then  by
Lemma~\ref{WeaklySymmetricDomesticStandardVsNonStandard}, $B$ is
necessarily nonstandard and $B=\Omega(n)$ for some $n$. But the
fact that $B$ is symmetric implies that  $char(K)=2$ which is a
contradiction.
\end{Proof}

In fact  we can extend (at least partially) the above result to
all weakly symmetric algebras of domestic representation type.

\begin{Thm}\label{domesticmaintheorem}
\begin{enumerate}
\item Two weakly symmetric  algebras of domestic
representation type are derived equivalent if and only if they are
stably equivalent of Morita type if and only if they are stably
equivalent.

\item  The class of   weakly symmetric
algebras  of standard domestic representation type  is closed
under stable equivalences, hence under stably equivalences of
Morita type, and derived equivalences.
\end{enumerate}
\end{Thm}

\begin{Proof} Since by
Lemma~\ref{WeaklySymmetricDomesticStandardVsNonStandard},
a standard weakly symmetric algebra of domestic type
cannot be stably equivalent to a nonstandard one,
Theorems~\ref{DomesticStandardSingular},
~\ref{DomesticStandardNonSingular} and ~\ref{DomesticNonStandard}
imply the first assertion.

For the second assertion, let $A$ be an indecomposable algebra
stably equivalent to a weakly symmetric standard algebra $B$ of
domestic representation type in the list of algebras in
Theorems~\ref{DomesticStandardSingular} and
~\ref{DomesticStandardNonSingular}. By \cite[Corollary 2 and the
discussions afterwards]{KrauseZwara}, $A$ is of domestic type.  As
in the proof of Proposition~\ref{ClassificationSymmetricDomestic},
$A$ is still self-injective. Since by
Proposition~\ref{ClassificationSymmetricDomestic}(2), the class of
symmetric algebras of domestic type is closed under stable
equivalences, one can assume that $B$ is not symmetric.  Now
$B=A(\lambda)$ with $\lambda\not\in\{ 0, 1\}$. Since
$B=A(\lambda)$ is a local special biserial algebra, so is $A$ by
\cite[Theorems 0.1 and 7.3]{Pogorzaly}. The algebra $A$ is thus
weakly symmetric.

 \end{Proof}

\begin{Rem}
\begin{enumerate}
\item If the class of  self-injective
nonstandard algebras of domestic type is closed under under stable
equivalences (hence under stably equivalences of Morita type  and
derived equivalences),  so is    the class of weakly symmetric
algebras of   domestic representation type.

\item The second assertion of the above theorem answers a question
of \cite[Page 47]{BocianHolmSkowronski2004}.

\item If we have a complete Morita equivalence classification of all
self-injective algebras of domestic type, the classification up to
stable equivalence (of Morita type) might be feasible as well.
\end{enumerate}
\end{Rem}

As a consequence,   the Auslander-Reiten conjecture holds for this
class of algebras.

\begin{Cor}
\begin{enumerate}
\item Let $A$ be an indecomposable algebra stably equivalent
to an indecomposable   symmetric algebra $B$ of
domestic type. Then $A$ has the same number of simple modules as
$B$.

\item Let $A$ be an indecomposable algebra stably equivalent
to an indecomposable weakly symmetric standard algebra $B$ of
domestic type. Then $A$ has the same number of simple modules as
$B$.
\end{enumerate}
\end{Cor}

\begin{Rem} \rm If the class of  self-injective
nonstandard algebras of domestic type is closed under under stable
equivalences, then Auslander-Reiten conjecture holds for a stable
equivalence between two weakly symmetric algebras of domestic
type.
\end{Rem}

\section{Non-domestic self-injective  algebras of polynomial growth}

\label{polysection}

The   derived equivalence classification of the standard (resp.
non-standard) non-domestic symmetric algebras of polynomial growth is
achieved in \cite[Page 653 Theorem]{BialkowskiHolmSkowronski2003a}
(resp.\cite[Theorem 3.1]{BialkowskiHolmSkowronski2003b}). They
give a complete list of representatives of derived equivalence
classes in terms of quiver with relations. We recall the derived
norm forms following \cite{HolmSkowronski2009}.

We have five quivers and various relations.

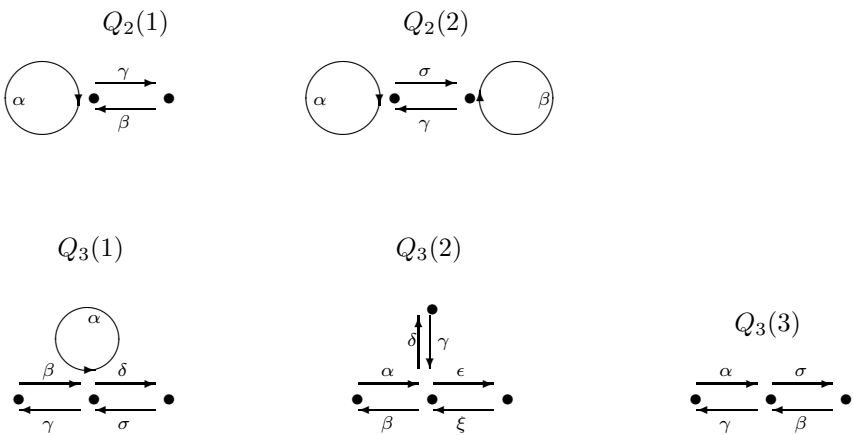
\unitlength1cm

\begin{picture}(15,6.5)

\put(1.2,6){$Q_2(1)$}
\put(1,5){$\bullet$}
\put(2,5){$\bullet$}
\put(1.1,5.3){\vector(1,0){.8}}
\put(1.9,4.95){\vector(-1,0){.8}}
\put(.4,5.1){\circle{1}}
\put(.88,5.1){\vector(0,-1){.1}}
\put(1.4,5.4){\scriptsize $\gamma$}
\put(1.4,4.7){\scriptsize $\beta$}
\put(0,5){\scriptsize $\alpha$}

\put(5.2,6){$Q_2(2)$}
\put(5,5){$\bullet$}
\put(6,5){$\bullet$}
\put(5.1,5.3){\vector(1,0){.8}}
\put(5.9,4.95){\vector(-1,0){.8}}
\put(4.4,5.1){\circle{1}}
\put(4.88,5.1){\vector(0,-1){.1}}
\put(6.7,5.1){\circle{1}}
\put(6.22,5.1){\vector(0,1){.1}}
\put(5.4,5.4){\scriptsize $\sigma$}
\put(5.4,4.7){\scriptsize $\gamma$}
\put(4,5){\scriptsize $\alpha$}
\put(6.99,5){\scriptsize $\beta$}

\put(9.6,2){$Q_3(3)$}
\put(9,1){$\bullet$}
\put(10,1){$\bullet$}
\put(11,1){$\bullet$}
\put(9.1,1.3){\vector(1,0){.8}}
\put(9.9,.95){\vector(-1,0){.8}}
\put(10.1,1.3){\vector(1,0){.8}}
\put(10.9,.95){\vector(-1,0){.8}}
\put(9.4,1.4){\scriptsize $\alpha$}
\put(9.4,.7){\scriptsize $\gamma$}
\put(10.4,1.4){\scriptsize $\sigma$}
\put(10.4,.7){\scriptsize $\beta$}

\put(.6,3){$Q_3(1)$}
\put(0,1){$\bullet$}
\put(1,1){$\bullet$}
\put(2,1){$\bullet$}
\put(.1,1.3){\vector(1,0){.8}}
\put(.9,.95){\vector(-1,0){.8}}
\put(1.1,1.3){\vector(1,0){.8}}
\put(1.9,.95){\vector(-1,0){.8}}
\put(1,1.9){\circle{.8}}
\put(1,1.48){\vector(1,0){.1}}
\put(.4,1.4){\scriptsize $\beta$}
\put(.4,.7){\scriptsize $\gamma$}
\put(1.4,1.4){\scriptsize $\delta$}
\put(1.4,.7){\scriptsize $\sigma$}
\put(1,2.1){\scriptsize $\alpha$}

\put(5.1,3){$Q_3(2)$}
\put(4.5,1){$\bullet$}
\put(5.5,1){$\bullet$}
\put(6.5,1){$\bullet$}
\put(4.6,1.3){\vector(1,0){.8}}
\put(5.4,.95){\vector(-1,0){.8}}
\put(5.6,1.3){\vector(1,0){.8}}
\put(6.4,.95){\vector(-1,0){.8}}
\put(5.5,2.2){$\bullet$}
\put(5.4,1.5){\vector(0,1){.7}}
\put(5.55,2.2){\vector(0,-1){.7}}
\put(4.9,1.4){\scriptsize $\alpha$}
\put(4.9,.7){\scriptsize $\beta$}
\put(5.9,1.4){\scriptsize $\epsilon$}
\put(5.9,.7){\scriptsize $\xi$}
\put(5.25,1.8){\scriptsize $\delta$}
\put(5.65,1.8){\scriptsize $\gamma$}

\end{picture}

Define for any $\lambda\in K\setminus\{0,1\}$ the algebras

\begin{eqnarray*}
\Lambda_2&:=&KQ_2(1)/(\alpha^2\gamma,\beta\alpha^2,\gamma\beta\gamma,
\beta\gamma\beta,\beta\gamma-\beta\alpha\gamma,\alpha^3-\gamma\beta)\\
\Lambda_2'&:=&KQ_2(1)/(\alpha^2\gamma,\beta\alpha^2,\beta\gamma,
\alpha^3-\gamma\beta)\\
\Lambda_3(\lambda)&:=&KQ_2(2)/(\alpha^4,\gamma\alpha^2,\alpha^2\sigma,
\alpha^2-\sigma\gamma-\alpha^3,\lambda\beta^2-\gamma\sigma,
\gamma\alpha-\beta\gamma,\sigma\beta-\alpha\sigma)\\
\Lambda_3'(\lambda)&:=&KQ_2(2)/(\alpha^2-\sigma\gamma,
\lambda\beta^2-\gamma\sigma,\gamma\alpha-\beta\gamma,
\sigma\beta-\alpha\sigma)\\
\\
\Lambda_5&:=&KQ_3(1)/(\alpha^2-\gamma\beta,\alpha^3-\delta\sigma,
\beta\delta,\sigma\gamma,\alpha\delta,\sigma\alpha,\gamma\beta\gamma,
\beta\gamma\beta,\beta\gamma-\beta\alpha\gamma)\\
\Lambda_5'&:=&KQ_3(1)/(\alpha^2-\gamma\beta,\alpha^3-\delta\sigma,
\beta\delta,\sigma\gamma,\alpha\delta,\sigma\alpha,\beta\gamma)\\
\Lambda_9&:=&KQ_3(2)/(\beta\alpha+\delta\gamma+\epsilon\xi,\gamma\delta,
\xi\epsilon,\alpha\beta\alpha,\beta\alpha\beta,
\alpha\beta-\alpha\delta\gamma\beta)\\
\Lambda_9'&:=&KQ_3(2)/(\beta\alpha+\delta\gamma+\epsilon\xi,\gamma\delta,
\xi\epsilon,\alpha\beta)\\
A_1(\lambda)&:=&KQ_3(3)/(\alpha\gamma\alpha-\alpha\sigma\beta,
\beta\gamma\alpha-\lambda\cdot\beta\sigma\beta,
\gamma\alpha\gamma-\sigma\beta\gamma,
\gamma\alpha\sigma-\lambda\sigma\beta\sigma)\\
A_4&:=&KQ_3(2)/(\beta\alpha+\delta\gamma+\epsilon\xi,
\alpha\beta,\gamma\epsilon,\xi\delta)
\end{eqnarray*}

Further, denote the trivial extensions of tubular canonical algebras as usual, in particular
$$\Lambda(2, 2, 2, 2, \lambda)=T(C(2,2,2,\lambda))\mbox{ for
}\lambda\in  K \setminus \{0, 1\}$$
$$\Lambda(3, 3, 3)=T(C(3,3,3))$$
$$\Lambda(2, 4, 4)=T(C(2,4,4))$$
$$\Lambda(2, 3, 6)=T(C(2,3,6)).$$

The precise quivers with relations of the trivial extensions
of tubular canonical algebras in question are displayed in
\cite{HolmSkowronski2009}.
We will refrain from presenting them here since we do not really need this
information in such details.

\subsection{Weakly symmetric standard non-domestic polynomial growth algebras}

\begin{Thm}
\cite[Page 653 Theorem]{BialkowskiHolmSkowronski2003a} Let $A$ be an
indecomposable standard non-domestic weakly symmetric algebra of
polynomial growth. Then $A$ is derived equivalent to one of the
following algebras:
\begin{itemize} \item  two
simple modules: $\Lambda_2'$ and  $\Lambda_3'(\lambda), \lambda
\in K \backslash\{0, 1\}$;

\item  three simple modules: $\Lambda_5'$ and $A_1(\lambda),
\lambda \in  K \backslash \{0, 1\}$;

\item four simple modules: $\Lambda_9'$  and $A_4$;

\item six simple modules: $\Lambda(2, 2, 2, 2, \lambda), \lambda
\in  K \backslash \{0, 1\}$;

\item eight simple modules: $\Lambda(3, 3, 3)$;

\item nine simple modules: $\Lambda(2, 4, 4)$;

\item ten simple modules: $\Lambda(2, 3, 6)$.
\end{itemize}
\end{Thm}

The above classification is complete up to the scalar problems in
$\Lambda_3'(\lambda)$ and in
$A_1(\lambda)$. 
Remark that except $\Lambda_9'$ in case the characteristic
of the base field $K$ is different from $2$, all
algebras are symmetric.

We shall prove

\begin{Prop}\label{NonDomesticStandard}
The classification of  indecomposable standard non-domestic
weakly symmetric algebras of polynomial growth up to stable
equivalences of Morita type coincide with the derived equivalence
classification, modulo the scalar problems in
$\Lambda_3'(\lambda), \lambda \in K \backslash\{0, 1\}$ and in
$A_1(\lambda), \lambda \in  K \backslash \{0, 1\}$.
\end{Prop}

\begin{Proof}
Since a derived equivalence between self-injective algebras
induces a stable equivalence of Morita type (\cite{KellerVossieck},
\cite{Rickard1989}), we only need to show that two algebras from
the above list which are not derived equivalent are not stably
equivalent of Morita type, either.

Since the property of being symmetric is invariant under a stable
equivalence of Morita type (\cite[Corollary 2.4]{Liu2008}), in
case the characteristic of the base field is different from $2$,
$\Lambda_9'$ cannot be stably equivalent of Morita
type to any of the other algebras. We can concentrate on
the remaining symmetric algebras in the list.

Now the algebras of at least six simple modules are trivial
extensions of canonical algebras of tubular type
$(2, 2, 2, 2)$,  $ (3, 3, 3)$, $(2, 4, 4)$,  $(2, 3, 6)$
respectively. They have singular Cartan matrices, while all
other algebras have non-singular Cartan matrix. Since, by
a result of Xi \cite[Proposition 5.1]{Xi2008} the
absolute value of the Cartan determinant  is preserved by a
stable equivalence of Morita type, we can consider separately
those algebras of at most four simple modules and those of at
least six simple modules.

For trivial extension cases, by \cite[Propositions 5.1 and
5.2]{BialkowskiHolmSkowronski2003a}, two algebras of this type are
stably equivalent if and only if they are derived equivalent. So
the number of simple modules, which is an invariant under derived
equivalence, distinguishes them.

Now suppose two indecomposable standard non-domestic
weakly symmetric algebras of polynomial growth are stably
equivalent of Morita type and have non-singular
Cartan matrices. Then the algebras are algebras in the above list.
The following table
gives some stable invariants of Morita type which distinguish two algebras
in the above list.
$$\begin{array}{c||c|c|c|c|c|c}
Algebra\  A & \Lambda_2' &\Lambda_3'(\lambda)& \Lambda_5'&
A_1(\lambda)& \Lambda_9'& A_4\\
\hline
det\  C_A& 6 & 12 & 6 & 16 & 4 & 12 \\
dim  (Z(A)/R(A))& 3 & 4 & 2 &2 & 1 & 2
\end{array}$$

\end{Proof}

\subsection{Selfinjective non-standard non-domestic polynomial growth algebras}

Now we turn to  non-domestic non-standard selfinjective algebras
of polynomial growth.  Recall the following

\begin{Thm} (Bia\l kowski, Holm, Skowro\'nski
\cite[Theorem 3.1]{BialkowskiHolmSkowronski2003b})
Let $A$ be an indecomposable non-standard non-domestic
self-injective  algebra of polynomial growth. Then we get the following.
\begin{itemize}
\item If the base field is of characteristic $3$, then $A$ is derived equivalent to $\Lambda_2$.
\item else the base field is of characteristic $2$ and then
\begin{itemize}
\item if $A$ has
two simple modules, $A$ is derived equivalent to $\Lambda_3(\lambda),
\lambda \in K \setminus\{0, 1\}$
\item  if $A$ has three simple modules, $A$ is derived equivalent to $\Lambda_5$
\item if $A$ has four simple modules, $A$ is derived equivalent to $\Lambda_9$

\item else $A$ has five simple modules and then $A$ is derived equivalent to $\Lambda_{10}$
\end{itemize}
\end{itemize}
\end{Thm}

If the base field is of characteristic $2$ and $\lambda\neq\lambda'$,
we do not know if $\Lambda_3(\lambda)$ is
derived equivalent to $\Lambda_3(\lambda')$ or not.
We call this question the scalar problem.

The above classification is complete up to the scalar problem
in case the base field is of characteristic $2$ for the algebra
$\Lambda_3(\lambda), \lambda \in K \backslash\{0, 1\}$.

Remark that
the algebra $\Lambda_{10}$ is not   symmetric (even not weakly
symmetric), so it is not stably equivalent of Morita type to any
other algebra in the above list.

We can now prove
\begin{Prop}\label{NonDomesticNonStandard}The classification  of  indecomposable
non-standard non-domestic self-injective algebras of polynomial
growth up to stable equivalences of Morita type coincides with the
derived equivalence classification, modulo the scalar problem in
$\Lambda_3(\lambda)$.

\end{Prop}

\begin{Proof} We only need to consider the case of characteristic
two, since if the characteristic of $K$ is $3$, there is only one algebra. But notice
$det\  C_{\Lambda_3(\lambda)}=12$,  $det\  C_{\Lambda_5}=6$ and
$det\  C_{\Lambda_9}=4$.
\end{Proof}

\subsection{Standard algebras versus non-standard algebras}

Now one can prove that a non-standard algebra cannot be stably
equivalent to a standard one. The case of a derived equivalence is
proved by Holm-Skowro\'nski in \cite{HolmSkowronski2009}.

\begin{Thm}\cite[Main Theorem]{HolmSkowronski2009}
Let $A$ be a standard self-injective algebra of polynomial growth
and let $\Lambda$  be an
indecomposable,  nonstandard, non-domestic, symmetric algebra of
polynomial growth. Then $A$ and $\Lambda$ are not derived
equivalent.
\end{Thm}

\begin{Prop}\label{NonDomesticStandardVsNonStandard}
Let $A$ be an indecomposable standard weakly symmetric algebra and
$\Lambda$ be an  indecomposable  nonstandard  non-domestic
self-injective algebra of polynomial growth. Then $A$ and
$\Lambda$ are not stably equivalent of Morita type.
\end{Prop}

\begin{Proof} Suppose that $A$ and $\Lambda$ are stably equivalent
of Morita type. By Krause \cite[Page 605 Corollary]{Krause}, $A$ is also
non-domestic of polynomial growth.

We only need to consider
symmetric algebras. In fact, the only non-symmetric nonstandard,
non-domestic   self-injective algebra of polynomial growth
$\Lambda_{10}$ is only defined when the characteristic of the base field is
$2$, but in case of
characteristic $2$, all standard non-domestic weakly symmetric
algebras of polynomial growth are symmetric, as $\Lambda_9'$ is
non-symmetric if and only if characteristic of the base field is different from
$2$.

Now we consider symmetric algebras. In fact, the method of the
proof of \cite[Main Theorem]{HolmSkowronski2009} works. Since two
algebras which are  stably equivalent of Morita type have
isomorphic stable Auslander-Reiten quivers, by the shape of the
stable Auslander-Reiten quivers (see the second paragraph of
Section 4 of Holm-Skowro\'nski \cite{HolmSkowronski2009}),
we proceed by remarking the following facts.
\begin{enumerate}
\item If $K$ is of characteristic $3$, the algebras $\Lambda_2$ and
$\Lambda_2'$ are not stably equivalent of Morita type.  Indeed
Holm-Skowro\'nski \cite{HolmSkowronski2009} have shown that
$$
dim\ Z(\Lambda_2')/T_1(\Lambda_2')^{\perp}=3\neq 2
=dim\ Z(\Lambda_2)/T_1(\Lambda_2)^{\perp}.
$$

\item If $K$ is of characteristic $2$ and
$\lambda,\mu \in  K \setminus \{0, 1\}$, the algebras $\Lambda_3(\lambda)$ and
$\Lambda_3(\mu)'$  are not stably equivalent of Morita type. Indeed
Holm-Skowro\'nski \cite{HolmSkowronski2009} have shown that
$$
dim\ Z(\Lambda_3'(\lambda))/T_1(\Lambda_3'(\lambda))^{\perp}=3\neq
2=dim\ Z(\Lambda_3(\mu))/T_1(\Lambda_3(\mu))^{\perp}.
$$

\item If $K$ is of characteristic $2$, the algebras $\Lambda_5$ and
$\Lambda_5'$ are not stably equivalent of Morita type. Indeed
Holm-Skowro\'nski \cite{HolmSkowronski2009} have shown that
$$dim
Z(\Lambda_5')/T_1(\Lambda_3')^{\perp}=0\neq dim
Z(\Lambda_5)/T_1(\Lambda_3)^{\perp}.$$

\item If $K$ is of characteristic $2$, the algebras $\Lambda_9$ and
$\Lambda_9'$ are not stably equivalent of Morita type.  Indeed
Holm-Skowro\'nski \cite{HolmSkowronski2009} have shown that
$$dim HH^2(\Lambda_9)=4\neq 3=dim HH^2(\Lambda_9').$$
\end{enumerate}

Since $Z(A)/T_1(A)^{\perp}$ and $HH^2(A)$ are invariants under
stable equivalences of Morita type, this completes the proof.
\end{Proof}

\subsection{Stable equivalence classification of weakly
symmetric non-domestic polynomial growth algebras}

Combining Propositions~\ref{NonDomesticStandard},
\ref{NonDomesticNonStandard} and
\ref{NonDomesticStandardVsNonStandard}, we have

\begin{Thm} \label{NonDomesticClassification}
The classification  of  indecomposable
non-domestic weakly symmetric algebras of polynomial growth
up to stable equivalences of Morita type coincides with the derived
equivalence classification, up to the above mentioned scalar problems.
\end{Thm}

As a consequence, we can prove a special case of Auslander-Reiten
conjecture.

\begin{Thm} \label{ARNonDomestic}
\begin{enumerate}
\item Let $A$ be an indecomposable
algebra which is stably equivalent of Morita type to
   an indecomposable
non-domestic  symmetric algebra  $\Lambda$   of polynomial growth.
Then $A$ and $\Lambda$ have the same number of simple modules.

\item Let $A$ and  $\Lambda$  be two indecomposable algebras which
are both  non-domestic  weakly symmetric algebra of polynomial
growth or which are both non-standard self-injective algebras of
polynomial growth. If they are stably equivalent of Morita type,
then $A$ and $\Lambda$ have the same number of simple modules.
\end{enumerate}
\end{Thm}

\begin{Proof} (1). By Krause \cite[Page 605 Corollary]{Krause}
and Liu \cite[Corollary 2.4]{Liu2008},
$A$ is also symmetric, non-domestic, and of polynomial growth.
Then one can apply Theorem~\ref{NonDomesticClassification}, by
noticing that the existing scalar problems all occur in families
of algebras for which different scalars yield algebras having
the same number of simple modules.

(2) is a direct consequence of
Theorem~\ref{NonDomesticClassification} and
Proposition~\ref{NonDomesticNonStandard}.

\end{Proof}

\subsection{Periodic algebras}

Recently, Karin Erdmann and Andrzej Skowro\'nski
\cite{ErdmannSkowronski2010} have shown that
a non-simple indecomposable symmetric algebra $A$ is tame with
periodic modules if and only if $A$ belongs to one of the following
classes of algebras:
a representation-finite symmetric algebra, a non domestic symmetric
algebra of polynomial growth, or an algebra of quaternion type (in
the sense of \cite{ErdmannLNM}). Note that these three classes of
algebras are closed under stable equivalences of Morita type, by
Krause \cite[Page 605 Corollary]{Krause} and Liu \cite[Corollary
2.4]{Liu2008}.

In \cite{Asashiba} Asashiba, in combination with Holm-Skowro\'nski \cite{HolmSkowronski2006}, the statement displayed in
\cite[Theorem 3.4, Theorem 3.5, Theorem 3.6]{ErdmannSkowronskiperiodic}
proved that for
the class of representation-finite 
selfinjective
indecomposable algebras are derived equivalent if and only if they are stable
equivalent. 
In our  recent paper
\cite[Theorem 7.1]{ZhouZimmermannTameBlocks},
we proved that for the class of algebras of quaternion
type, derived equivalence classification coincide with the
classification up to stable equivalences of Morita type
(up to some scalar problems).
Therefore, combining the above Theorem~\ref{ARNonDomestic}, we
showed the following

\begin{Thm}
\begin{enumerate}
\item
    The class of  tame symmetric algebras  with
periodic modules is closed under stable equivalences of Morita
type, in particular under derived equivalences;

\item
For tame symmetric algebras  with periodic modules
derived equivalence classification coincide with the
classification up to stable equivalences of Morita type (up to
some scalar problems);

\item   The Auslander-Reiten conjecture holds for a stable
equivalence of Morita type between two tame symmetric algebra with
periodic modules.
\end{enumerate}
\end{Thm}

\section{Concluding remarks}

In Section~\ref{polysection}
we classified symmetric algebras of polynomial growth
only up to
stable equivalences of Morita type, whereas the results in Section~\ref{weaklysymmdomesticsection} concern stable equivalences in general.
It would be most interesting to get a stable equivalence classification
for general selfinjective algebras of polynomial growth,
in particular the validity of the Auslander-Reiten conjecture in general.
The crucial point for this more general statement would be
to first classify algebras which are stably equivalent (of Morita type) to
$\Lambda_9'$ in case the characteristic of $K$ is different from $2$
or to $\Lambda_{10}$, since these are the selfinjective polynomial growth
algebras which are  not symmetric. The second step is
to classify  weakly symmetric algebras of polynomial growth
up to stable equivalences.

If the for the class of selfinjective algebras
of polynomial growth the derived equivalence classification would coincide with the
classification up to stable equivalences (of Morita type) the Auslander-Reiten
conjecture for this class of algebras should follow.

Another point concerns the work of Martinez-Villa (\cite{Martinez-Villa}) of the
reduction of the Auslander-Reiten conjecture to selfinjective algebras.
If we can prove
that this reduction preserves representation type:
domestic type, polynomial growth etc,
then the conjectural result of the preceding paragraph would imply
the validity of the
Auslander-Reiten conjecture for algebras of polynomial growth.

The results of this paper present a second occurrence of a phenomenon of the same
kind: In our previous paper
\cite{ZhouZimmermannTameBlocks} we also obtained that the classification of a class of
tame symmetric algebras up to stable equivalence coincides with the derived equivalence
classification.
One might ask if two indecomposable tame symmetric algebras are stably equivalent of
Morita type if and only if they are derived equivalent.

\end{document}